# Mathematics in art, for art and as art

**Maria J. Esteban**

**Affiliation**: CNRS UMR 7534, PSL University, Université Paris-Dauphine, Place de Lattre de Tassigny, 75775 Paris 16, France
esteban@ceremade.dauphine.fr

**Abstract:** The fundamental role of mathematics as an inspiration for artists, but also as a tool for art creation, is presented in this paper following different art fields, like architecture, sculpture, painting, photography, literature and poetry, movie making and music. The historical viewpoint is completed with recent applications of mathematics to create art in the digital era. Finally, the article contains a discussion about the possibility of the mathematical creation being considered artistic.

## Introduction

In ancient times, and until three centuries ago, scholars were men and women who were learned in various scientific fields, they knew about mathematics, mechanics and astronomy and often were called mathematicians or even philosophers. At that time, science was not divided into different areas of specialization. Some of those men and women knew more of this or that, but their knowledge embraced many aspects of what today we would call science and technology. Mathematicians were asked by kings, emperors, and pharaohs to help designing and building palaces, fortresses, temples, imposing and awesome graves and whole cities. Their knowledge was important not only because of practical reasons, but also because of aesthetical ones. Mathematics, and especially geometry, helped in designing buildings that were solid, that would not collapse despite their impressive dimensions, but which were also beautiful, attractive to the eye and surprising by their huge dimensions and nevertheless, elegant and light looking. In ancient times, mathematics was born out of the need to solve concrete problems: how to build, how to keep accounts, how to measure. And those different needs initiated different parts of mathematics, like algebra, geometry, and later, calculus. Geometry was key in the architectural use of mathematics because it is concerned with shapes, but also with proportions, which are both of the utmost importance to build beautifully. But the solidity of buildings required also careful calculations that were often in the hands of mathematicians who worked as architects.

Saying that artists are men and women is not perfectly accurate. Indeed, Nature and animals can also be "artists" and "use" geometry in their designs. This is due to the principle of optimization, when one (animals, humans, Nature) aims to use a minimal amount of energy, building material, time or to spend a minimal amount of money, or to be as resistant under deformation as possible. Yes, when bees make their *honeycombs* (Fig. 1) using hexagonal patterns, *they are trying to minimize* the amount of wax they need. Fascinating, is not it? And more fascinating indeed is that hexagons are also found on bees' eyes. Actually, an open question in biology is to understand whether there is a link between the two phenomena or not. But hexagons are not found only in bees' honeycomb patterns, they are ubiquitous in the vegetal and mineral worlds too: as a single example, let me just speak of the incredible hexagonal patterns found in *basalt columns*, or *basalt organs* (Fig. 2), as the French call them. Nowadays, different hexagonal honeycomb materials are used in the design of carbon based very thin materials, like in graphene, and this is due to the need to get the best effective and resistant materials. Explicitly or implicitly, geometry and calculus are behind all those optimization problems.

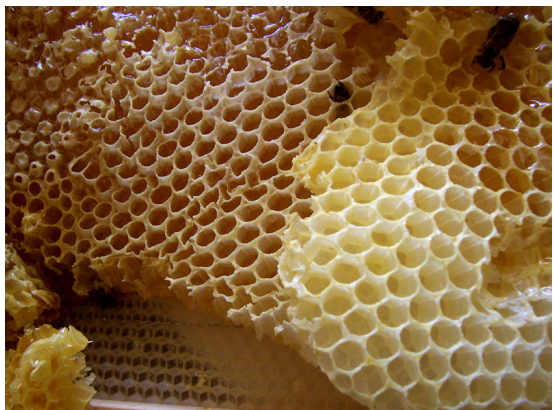
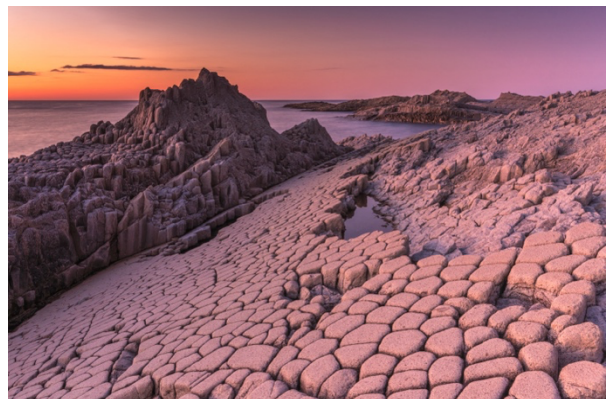

(1) Photo by Merdal / CC BY-SA 3.0          (2) Photo by Екатерина Васягина / CC BY-SA 4.0

Beyond architecture, the importance of mathematics in the art of painting appears much later, mainly when in the 15th century the Dutch artists wanted to paint landscapes looking realistic. For that purpose, they needed perspective. In the art of painting maybe the most important time when mathematics played a crucial role was when the notion of *linear perspective* was defined and systematized to represent three-dimensional landscapes on two-dimensional surfaces. This geometrical invention is often attributed to Brunelleschi, famous architect of the 15th century. Before him, other artists had understood that lines which converge to a point could create some impression of depth in paintings, like the Italian painters Duccio or Giotto in the 13th-14th centuries, but their work was more based on intuition than on a rigorous system. Chinese painters and Greek artists used other techniques to give an impression of depth and volume, but none of their methods was systematized and based on precise mathematical rules. Later, systematized perspective methods were tools in the hands of architects, sculptors and painters. Many centuries later, mathematics entered the artistic realm in a much more explicit way: in the last 100 years, many paintings are full of mathematical objects and architects and sculptors also use them widely in their designs.  For instance, mathematics was a source of inspiration and a tool in the development of *cubism*, in the *abstract art* movement, for the artists of the Bauhaus school and among the suprematists in the Soviet Union.  Later, among others, fractals can be found in the works of many artists and Vasarely's optical creations are also inspired by mathematics.

Another artistic field where mathematics is omnipresent is music, and this is also a very old phenomenon. Indeed, (musical) sound can be described using mathematics, and more concretely, the notions of *numbers*, *ratios* and *scales*. In ancient Greece harmony was considered a fundamental branch of physics, now known as *musical acoustics*. And similar approaches are found in old India and China. Therefore, mathematical fields like algebra and number theory have played an important role in musical composition. More recently, even fields like topology and differential geometry have been used in the analysis of music and its structure. For instance, the famous contemporary composer I. Xenakis has based some of his compositions on Markov chains and other mathematical objects and tools, as he describes it in [1]. Lately, even more modern algorithmic technologies are used to create *digital music* with the help of computers and a minimalist composer like Tom Johnson uses all kinds of mathematical concepts, results and theories as a tool for musical composition.

How not think of geometry when one looks at old Arabic tiles, old rugs and stucco decoration in Islamic architecture? Also, symmetry in shapes is one of the most inspiring impacts of mathematics in art, because symmetry is related to beauty and so, symmetric designs are omnipresent in art, in every culture, in every century.

All the above is about how mathematics can help and has helped in creating art. But there are two other aspects that I would like to discuss in this article. The first one is about the fact that mathematics can help with art not only creating it, but also studying, analyzing and repairing it. For instance, the *level set method* and other techniques based on the use of *nonlinear differential equations* are key in the analysis of images, to repair old videos, pictures, and sound recordings. Also, based on the *theory of wavelets*, the JPEG 2000 and the JPEG XS protocols are nowadays the main tools to compress pictures and all kinds of images and videos in an incredibly efficient manner.

Furthermore, mathematics-based methodologies and technologies are heavily involved in the art of making videos, movies and photography. Not to speak of how mathematical algorithms are behind the creation of many film effects in the movie industry, especially in animated movies, but not only. It can also be mentioned that mathematics has been used as a tool to analyze problems of authorship in different arts.

The second and final additional topic that I would like to discuss here is very far from the above, but I find it of importance: is mathematics itself an art? The creativity needed to do new mathematics is assimilable to the creativity of an artist in the classical definition of it?

This is a topic that I would like to discuss a bit together with the other aspects of creativity relying on mathematics.

The paper is organized as follows. In Section 1 we will discuss in more detail the relations between mathematics and music, with an important emphasis on how mathematics is changing the face of music today and how simultaneous education in both fields is becoming quite common in many universities. In Sections 2, 3 and 4 we will see in more detail the impact and importance of mathematics in architecture, painting and sculpture. Section 5 will explore the use of mathematics in the world of photography and the movie industry. After quickly discussing the link between mathematics, literature and poetry in section 6, in section 7 we will explore in which measure mathematical creativity can be considered an art by itself. Some conclusions can be found at the end the paper.

**1. Mathematics and music**

Since ancient times, music has been considered a mathematical art, in the sense that music's description is based on *scales* and *ratios* between numbers. In the Western tradition, Pythagoras was certainly the first to be interested in this relation. Many works, articles and books deal in great detail with examples coming from different centuries and regions of the World. They not only describe historic examples of the relation between mathematics and music, but also the technical reasons behind this strong link.

Interestingly, many mathematicians are excellent musicians and even enjoy a (semi-)professional status as music players. Music is certainly the artistic branch which is more spread among mathematicians. But what is interesting is that very well-known musicians also studied mathematics. See for instance https://en.wikipedia.org/wiki/Music_and_mathematics.

**1.1. Music composition and computational music**

Mathematics has always been present in music composition because it is the structure that each composer chooses for his/her compositions that is behind his/her creative spur. Without structure, without rhythm, music would be only noise. And of course, that structure is more or less explicitly based on mathematical rules. Numbers and structures offer the composer the building blocks on which he/she will act to finish his/her composition. As Sarah Hart says in her lecture on the mathematics of musical composition in Gresham College [2], "Music involves pitch, melody, and rhythm. The earliest music, before instruments were made that could make specific sounds, must have been made largely with rhythm: clapping, stamping, or repetitive vocalisations. Percussion instruments like drums, gongs, bells, exist in various forms in almost every culture. Later, along with our own voices, instruments were invented that could play defined pitches. Melody and harmony were developed. Music started to acquire forms, and to be written down. These forms developed along different lines in different cultures". It is well-known that Bach was able to write music in every key so successfully because mathematicians found better ways to calculate the 12th root of the number 2. This is related to the musical problem of dividing the octave into twelve equal intervals, which involves splitting sound waves into *ratios* rather than equal lengths.

As the mathematician and concert pianist E. Cheung said in her interview to the *uchicago news journal* [3], "During the Baroque period, a mathematical breakthrough inspired one of Cheng's favorite composers, Johann Sebastian Bach, to write *The Well-Tempered Clavier* (1722), his book of preludes and fugues in all 24 major and minor keys. That's why music before the Baroque time didn't really modulate. It always stayed in the same key. Because of the way that they tuned keyboards, if they moved a key, it would have sounded terrible".

In the twentieth century we could see a new phenomenon arising, because many composers consciously started using mathematical concepts and theories as bases for their creative processes. For instance, first rate composers like Olivier Messiaen, Iannis Xenakis, Franco Evangelisti, Pierre Boulez, Arvo Pärt, Steve Reich and Philip Glass have followed this path and in [4] we can find a very complete presentation of this trend.

The work of many contemporary musicians shows the power of mathematics to open new possibilities in music. Modern experiments with *digital computer music* are just the most recent example: in a computer one can easily change sounds to produce new music.
In the famous *Institut de recherche et coordination acoustique/musique* (IRCAM) in Paris, one can find research group topics like *Modèles algébriques, topologiques et catégoriels en musicologie computationnelle*, very surprising at first sight. You find also professor chairs and courses an PhDs with mathematics-related titles.

## 1.3    Math & Music education and research

Music is certainly based on a mathematical structure and so it is important for music students to understand enough mathematics so that they can appreciate better the structure of music. Many universities all over the world are offering joint degrees, or at least courses, in mathematics and music. And these courses do not only teach subjects in classical fields of mathematics, but also, they propose the students to get trained in computer programming which will allow them to analyze music and to create it. Since from the vibration of a guitar string and the analysis of a rhythmic pattern to the use of randomness in experimental music mathematics continues to impact music at every level, this kind of educational programs propose to teach students topics both in pure and applied mathematics. But even more, nowadays there are high schools with integrated courses in mathematics and music. And of course, machine learning and AI have also been put to contribution in what concerns music, and at IRCAM there are several research groups working on this topic. For instance, there are PhD theses linking mathematics and music, as the thesis of Paul Lascabettes on *Models for the Discovery of Musical Patterns and Structures, and for Performances Analysis*.

In the contemporary music community, a lot of activity is going on around the interaction of mathematics and music and there are also international journals concerned with this topic, like the *Journal of Mathematics and Music*, published by Francis and Taylor. Recently, in February 2024, a big conference took place in Singapore about *Mathemusical Encounters in Singapore: a Diderot Legacy.* But this kind of conferences have been going on for some time already. For instance, in 2011 took place in the *Centre Pompidou*, in Paris, the international *Mathematics and Computation in Music* Conference. And in 2022 took place in London the *IMA Conference on Mathematics in Music*. Exhibitions are also organized on this topic, like the *La La Lab – the Mathematics of Music*", which took place in Heidelberg in 2020.

There is a very large number of publications on this topic, like for instance, the *Mathemusical conversations: mathematics and computation in music performance and composition*, whose editors are J.B.L. Smith (Japan), E. Chew (UK) and G. Assayag (France) [5].

**2. Mathematics and architecture**

Since the antiquity, mathematics has been very important in architecture and for architects, and this for two reasons: to build beautiful, harmonious, and well-proportioned buildings, but also to build them in a way that they would stand up and would not collapse. When drawing the plans, when executing and finalizing them, architects, and the engineers working with and for them, use mathematics intentionally or unintentionally. But note that even if mathematics' role is fundamental in their work, very few will understand it. Even in buildings which look totally ageometric, like for instance, Gehry's main buildings, the engineers who have to interpret his drawings, have to make sure that the structures will stand firmly, and for that, numerous structural studies have to be made. Of course, they will use computer-based algorithms to manipulate the shapes and conceive stable structures, but those CAO algorithms are built themselves on equations and calculations containing a good amount of mathematics and mechanics.

When building a new structure, the main steps where mathematics plays an essential role are found in the structural calculations, in the geometry and spatial forms, in the search for beauty through symmetries and proportions; and finally, in the interior decoration of surfaces (think for instance about the tiles and the stucco decoration in Islamic architecture).

In ancient Greece, India, Mesopotamia, and Egypt, special proportions were used when building temples, mausoleums, pyramids, etc. and this often for religious reasons. In some cases, circles were important like in China and in some Buddhist buildings too. In India one can find constructions that remind us of fractal shapes.

In mathematics it is said that two quantities are in the *golden ratio* if their ratio is the same as the ratio of their sum to the larger of the two quantities. In Greece, and in particular in Euclidean geometry, the *golden ratio* was called *extreme or mean ratio*. And this proportion was considered to be so related to beauty that it was even designated as the *divine proportion* by Luca Pacioli, a contemporary of Leonardo da Vinci, and one of the greatest experts in accounting in the 16th century. The *golden ratio* was thought to be used to build the pyramids in Egypt and the Parthenon in Greece, but recent studies have shown that this was probably not the case, because this *ratio* was known around 300 BC, but probably not before that. In any case, proportions were very important for architects when designing the pyramids or the beautiful old Greek temples. Mathematics was also key when, building the Middle Age churches and cathedrals, especially when building the high and elegant ones of the Gothic period. And when the architects became too ambitious and went beyond what was reasonable from the stability point of view concerning arches, columns and roofs, some cathedrals collapsed, like for instance the Beauvais cathedral near Paris. In the Renaissance, Brunelleschi was a key figure with his discovery of the principles of *linear perspective*. Many buildings were erected following these principles. A large number of mathematicians made contributions to architecture at that time, and Leone Alberti wrote important texts on the matter. Leonardo da Vinci, who loved mathematics and learned a lot from Alberti, and who is not that well known as an architect, was not only a fabulous painter, but also an architect and an engineer, and because of these skills, he worked for the Duke of Milano and the King of France.

In the 17th century, an architect called C. Wren wrote that the geometric figures are naturally much more beautiful than the irregular ones. And much later, Le Corbusier said more or less the same: for him cubes, cylinders, spheres and pyramids were the best shapes to build beautifully. This, of course, has not been shared by many modern architects who decided to design their buildings using more complex geometries. Also, in modern times architects started using more sophisticated designs based on elaborated mathematical curves, like hyperbolic paraboloids and hyperboloids of revolution, catenoids, helicoids, and ruled surfaces in the case of Gaudi. Among the first structures built using hyperboloids is the Shukov tower in Russia (1896). More recently, a beautiful example is the Brasilia cathedral built by Niemeyer (1970).

Dome architecture started probably in ancient Rome, with the Pantheon in Rome and continued to develop in cathedrals built during the Renaissance, mainly in Italy. Later, domes are omnipresent in France, the USA, in Eastern European church architecture (as *onion domes*), and elsewhere. Geodesic domes are also ubiquitous nowadays. You can find them in the Eden project in the UK, or La Geode in Paris, or the National Center for the Arts in Beijing (Fig. 3), the Biosphere of Montreal, the Telus Sphere in Vancouver, the Matrimandir in India (Fig. 4), etc, etc.

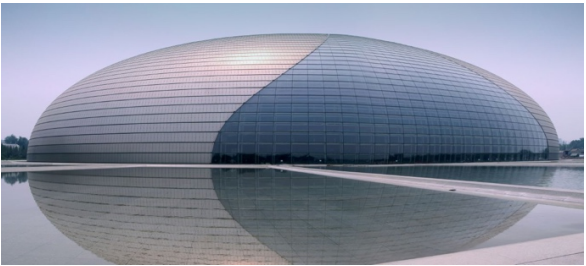 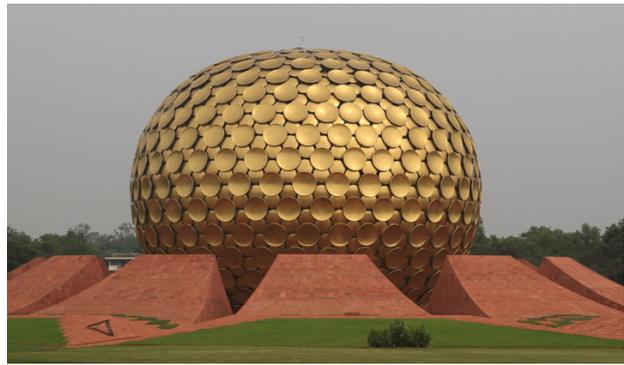

(3) Photo by Hui Lan / CC BY 2.0  (4) Photo by Ve.Balamurali / CC BY-SA 3.0, cropped

What is commonly called modern architecture spans mainly from 1920 to the 1980s. The main principles of that style were to construct simply looking and functional buildings, without unnecessary ornaments and using new technologies of constructions, and concrete, steel and glass as basic materials. Very famous examples of this style are the buildings of Niemeyer in Brasilia, Le Corbusier's houses, many buildings belonging to the Art Deco style, the iconic Fagus factory, built in Alfeld by A. Meyer and W. Gropius, the founder of the *Bauhaus*. Let us emphasize that modern architecture had been prepared mainly by American architects, for instance those belonging to the Prairie school, and among them, the most famous, F. Lloyd Wright. They had started building functionally, with simple designs, abandoning curves and using many straight lines.

The *Bauhaus*, which was created by Walter Gropius in 1919, and was linked to the *De Stijl* movement, synthetized very well the principles of simplicity and functionality. Among the participants in this group, we find not only architects, but also interior designers, furniture makers, photographers and, as we shall see later, many painters and sculptors. Their designs were very geometrical looking, simple, with many straight lines and rectangles as the main elements of construction. Another very visible architect belonging to the school of *rational modernism* and *international style* was Alvar Aalto.

The International exhibition of decorative arts, which took place in Paris in 1925, was an excellent showcase for this modern, simple, and geometrical architectural style, with linear, simple and white buildings.

Later, many buildings have been designed borrowing from geometry. It is impossible to describe or discuss all of them here. Very interesting modern designs related to mathematics include the Opera House in Sidney (Fig. 5) and the Lotus temple in India (Fig. 6), a lotus-shaped monument conceived with a combination of various geometrical shapes. It was designed by architect F. Sahba and the lotus flower had to be simplified into geometrical forms such as spheres, toroids, cones, and cylinders. Finally, impressive buildings with very surprising and beautiful geometrical designs have been also built by Zaha Hadid, who was influenced by the *suprematists* and the *constructivists* active in the Russian art scene almost one century before.

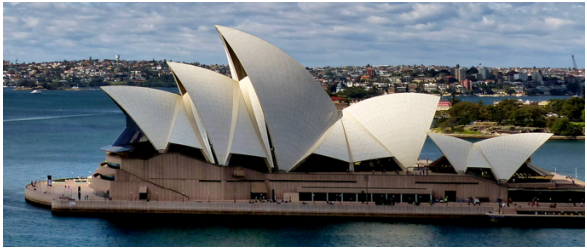 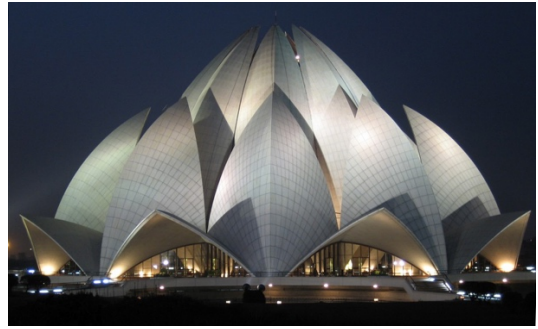

(5) [Photo](...) by [Bernard Spragg. NZ](...) / [CC0](...), cropped.   (6) [Photo](...) by [Vandelizer](...) / [CC BY 2.0](...)

Let us end this section speaking a bit about *tessellation* and *tiling*s. *Tilings* have been present in interior decoration from the very oldest times. The repetition of some periodic shapes created beautiful designs for walls and floors. They were used already by the Sumerians. We find this kind of periodic designs in most Islamic buildings, in tiles and in the stucco decoration of the walls and ceilings. They were also often used in Roman mosaics. Nowadays we find also aperiodic tilings for instance in the external walls of the famous National Aquatics Center in Beijing (Fig. 7). In this case, the tiling is not flat, it is composed of bubble-like elements. The external walls of this impressive building form a natural pattern of bubbles in soap lather. One by one, they are solutions of the *minimal Surface Plateau problem*, and the realization of this building was the first implementation of the solution to that geometric problem stated in 1887, which was solved only in 1993.

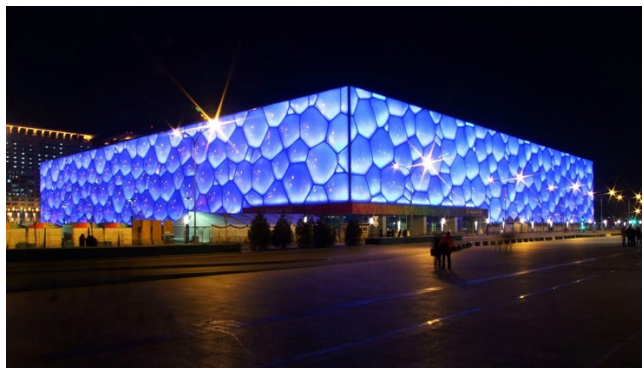

(7) [Photo](...) by [Charlie fong](...) / [Public Domain](...)

### 3. Mathematics and sculpture

The first formalized and mathematical approach about how to sculpt beautifully is Polykleitos' *Canon*, which described how to make elegant human sculptures. Even if the *Canon* has been lost, the rules described in it are widely well-known. Polykleitos' model was a nude male body, and he described its optimal proportions to ensure perfection. It is believed that his main example was his *Doryphoros* [6]. The main notions present in the *Canon* were symmetry, rhythm, and balance. An interesting feature was that these notions applied not only to the whole body, but to all its parts, one by one. The influence of the *Canon* on Greek, Roman and Renaissance sculpture has been enormous, facilitating the existence of thousands of very beautiful sculptures for many centuries.

In the 20th century, mathematics has influenced sculpture in several ways. The presence of mathematics in modern sculptures is more explicit, since many of them follow or represent geometrical objects and curves, some of them explicitly, some less so. See for instance an artist like Hiroshi Sugimoto that I came to know in the exhibition *Mathématiques, un dépaysement soudain* [7] presented in 2011 at the *Fondation Cartier pour l'Art Contemporain* in Paris. In that exhibition we could see an impressive sculpture of his representing a surface of revolution with constant negative curvature. A series of his works is actually called "Mathematical models" [8].

Nowadays, many artists use mathematics in their work, making sculptures that are often called *mathematical sculptures*, and which represent mathematical (often geometrical) objects. One can see them in museums, exhibitions and in the streets. One of these mathematical sculptures that I always liked because it is beautiful and because it was designed by mathematicians and made by some artist following their design, is the *Dodekaederstern* (Fig. 8), standing in front of the main entrance of the Mathematical Department of the University of Vienna. This geometric figure is defined in all its details by just one single algebraic equation, and it illustrates the interplay of various fields of mathematics such as invariant theory, algebraic geometry, and singularity theory.

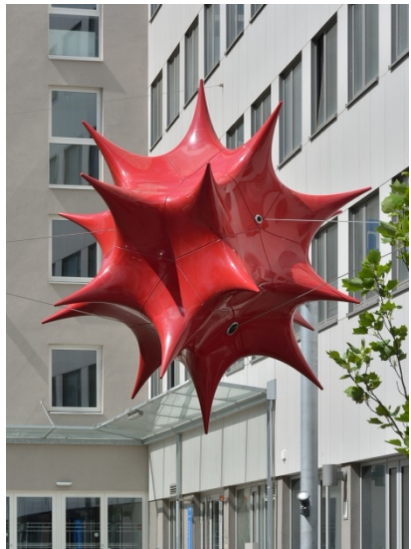

(8) Photo by Peter Haas / CC BY-SA 3.0

The *Infinite staircase* of Olafur Eliasson (2004) (Fig.9) sitting in the atrium of an office building in Munich, represents a continuous loop in the form of a double helix. The heights of the steps vary lightly to compensate for the curvature of the staircases. It is an elegant and surprising sculpture which required precise engineering to be able to stand in balance. In some sense, when one looks at this sculpture, one is reminded of the double helix staircase of the 16th Century Chateau de Chambord in France (Fig.10). This very original and famous staircase is attributed to Leonardo da Vinci.

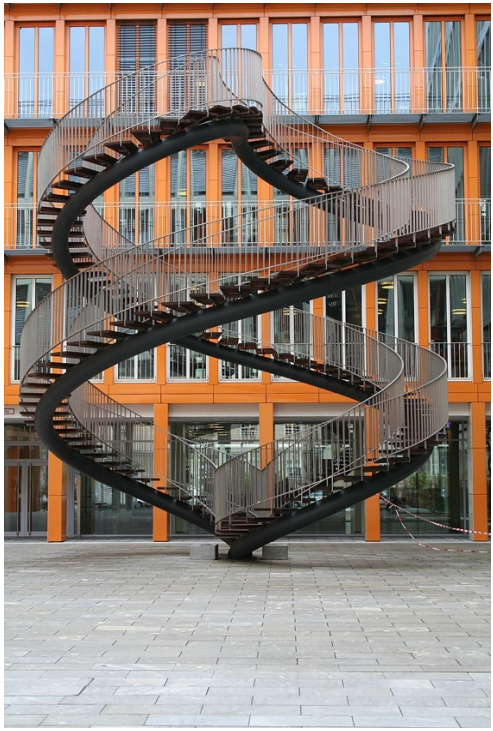 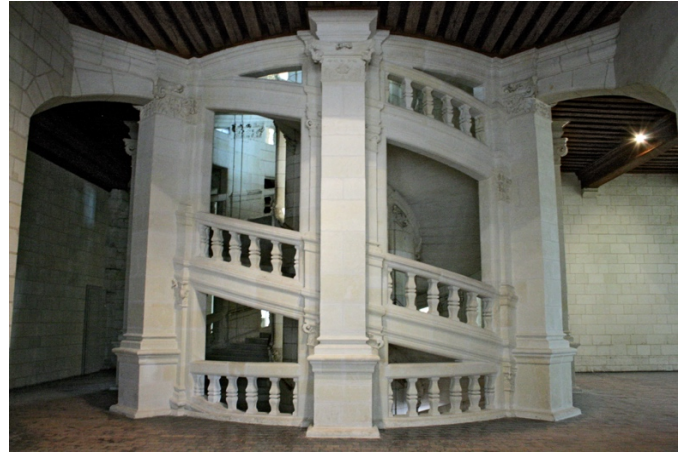

(9) [Photo](#) by [Oliver Raupach](#) / [CC BY-SA 2.5 DEED](#)

(10) [Photo](#) by Hélène Rival / [CC BY-SA 3.0](#)

*Mathematical sculptures* are ubiquitous in museums, galleries and in the streets of our modern cities, and the European Society for Mathematics and the Arts (ESMA) periodically announces exhibitions and activities where *mathematical sculptures* play a central role. But *mathematical sculptures* are not only made by artists. As an anecdote, let us say for instance that the game called *Mathematical origami* is played by children, youngsters, and adults, with contests organized to select the most intricate and beautiful of those objects. Also, many geometric objects and objects which are solutions of mathematical equations are leading to the creation of sculptures by using computers and 3D printers.

## 4. Mathematics in painting and design

### 4.1 Perspective in painting

In very old paintings and frescoes we can already acknowledge the efforts of their creators to give the objects a sense of distance from the observer, buildings and landscapes in the paintings by using various perspective tricks. In old Chinese and Indian paintings, perspective effects were already in place, but they were not systematized, not really based on precise mathematical rules. They are mostly based on the so-called *Oblique projection*, a simple type of technical drawing of graphical projection to represent three-dimensional objects on a flat two-dimensional image. The objects are not in perspective and so they do not correspond to any real view that one could see, but the technique yields somewhat convincing results.

Medieval artists in Europe, like those in the Islamic world and China, were aware of the general principle of varying the relative size of elements according to distance, but without being based on a systematic formalized theory. Byzantine art used the so-called *reverse perspective* convention for the setting of principal figures. In *reverse perspective* the objects depicted in a scene are placed between the projective point and the viewing plane. Objects farther away from the viewing plane are drawn as larger, and closer objects are drawn as smaller.

The first systematized approach to perspective is the so-called *linear perspective*, a rigorous geometrical theory probably initiated by the architect Filippo Brunelleschi in the Italian Renaissance, in the 15$^{th}$ Century. In *linear perspective*, all lines appear to converge toward a vanishing point in the painting (sometimes there can be multiple vanishing points) and the main effect of *linear perspective* is that objects appear smaller as their distance from the observer increases. Brunelleschi understood this technique and ran a series of experiments to explain the interest of *linear perspective*, but he did not publish the mathematics behind it. Piero della Francesca elaborated on *De pictura* in his *De Prospectiva pingendi in* the 1470s, making many references to Euclid. The main theoretical work in the subject is Luca Pacioli's *Divina Proporzione* (1509) [9], illustrated by Leonardo da Vinci, which summarizes the use of perspective in painting, including much of Della Francesca's treatise. Leonardo applied *one-point perspective* to some of his works.

*Linear perspective* is an advance with respect to the previous non-formalized methods but has its limitations. Distortions appear, and for instance, a sphere drawn in *linear perspective* will be stretched into an ellipse. These apparent distortions are more pronounced away from the center of the image and artists sometimes tended to manually correct those distortions, thus abandoning the formal method.

*Two-point perspective* was demonstrated as early as 1525 by Albrecht Dürer, who studied perspective by reading della Francesca's and Pacioli's works. The *three-point perspective* uses three vanishing points to portray objects or scenes with notable depth or height. This technique is commonly used when observing tall structures like skyscrapers or deep canyons, where a third vanishing point appears above or below the horizon line.

*Curvilinear perspective*, also called *five-point perspective*, is another graphical projection used to represent three-dimensional objects on two-dimensional surfaces. It was formally codified in 1968 by the artists and art historians André Barre and Albert Flocon in the book *La Perspective curviligne*. It seems that the intuition behind this work comes from the book *Grafiek en tekeningen*, by M.C. Escher. This kind of perspective is sometimes colloquially called *fisheye perspective*, by analogy to a fisheye mirror. But *five-point perspective* was used much earlier in an approximate manner. This was for instance the case of the Arnolfini portrait of the Flemish painter Van Eyck.

Summarizing, perspective in art is a technique that communicates depth, space, and realism on a two-dimensional medium and that is based on geometric principles.

## 4.2 Geometry in painting

Piet Mondrian was interested in theosophy, mathematics and geometry. His contacts with theosophy pushed him towards simplicity in shapes and colors, and therefore, towards abstraction. Around 1911 he saw some paintings of Braque belonging to the *analytic cubism* period. This led him to start decomposing landscapes in that style, and already in 1912, he painted destructured landscapes very close to Braques' style in grey and ochre colors (Fig. 11). In the following years he found a new pictural language, and he became one of the pioneers of abstract art, together with Kandinsky, Kupka, Léger, Picabia and Robert and Sonia Delaunay. These painters chose to paint in simple, pure colors, and using geometric objects, like lines, circles, squares and rectangles among others. Mondrian went to big extremes, painting very simple sets of lines which cross orthogonally. His colors became also pure and striking (Fig. 12).

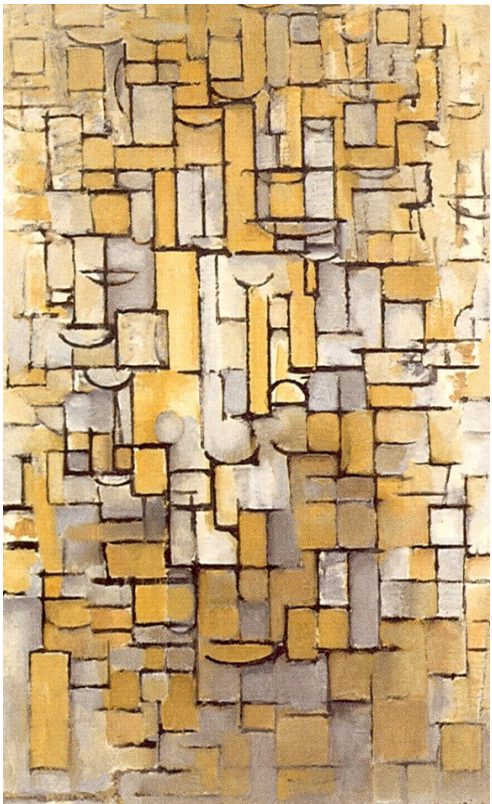
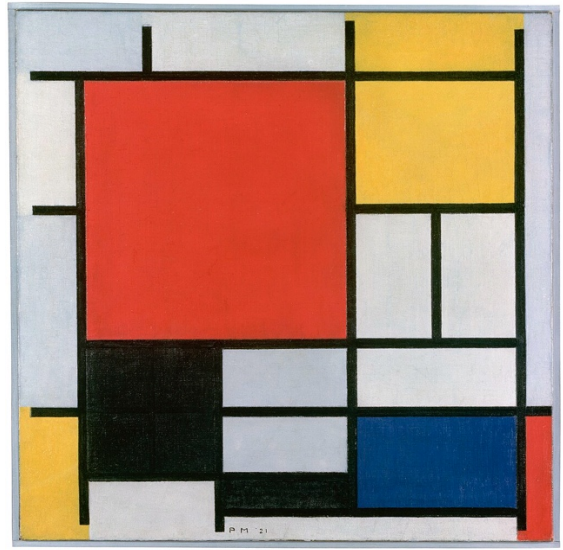

(11) Photo /Domaine public    (12) Photo /Domaine public

Mondrian was part of the *De Stijl* Dutch group, together with Theo van Doesburg and others. Note that in 1916 van Doesburg painted a geometric painting called *Composition I (Still Life)*, following the same principles as Mondrian.

As mentioned previously, in the years 1919 to 1933, the *Bauhaus* was a group of artists willing to act in all forms of art, architecture, painting, sculpture, interior design, pottery, etc. The main features identified in *Bauhaus* style are simple geometric shapes like rectangles and spheres, without elaborate decorations. Kandinsky (Fig. 13), also a theosophist, tending more and more towards simplicity, was invited to join the *Bauhaus* in 1922, and there he was in charge of painting classes and workshops, where he taught his theory of colors and shapes. Geometrical elements took on increasing importance in both his teaching and painting—particularly the circle, half-circle, the angle, straight lines and curves. In 1926 he published a theoretical book entitled *Point and Line to Plane*.

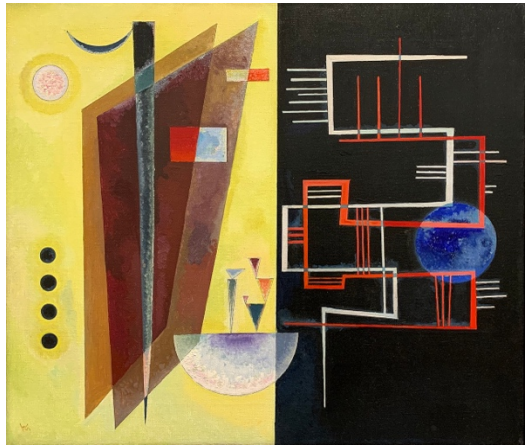

(13) Photo by Benoît Prieur et un auteur supplémentaire / Domaine public

In 1915-16, also pushed by philosophical and spiritual ideas, K. Malevich (Fig. 14) laid down the foundations of *Suprematism,* an art movement focused on the fundamentals of geometry (circles, squares, rectangles), and he painted in a limited range of colors. The influence of this current has lived for very long, since in the 21$^{st}$ century, the study of the architect Zaha Hadid, that we mentioned in the section about math and architecture, studied the works and theories around Malevich.

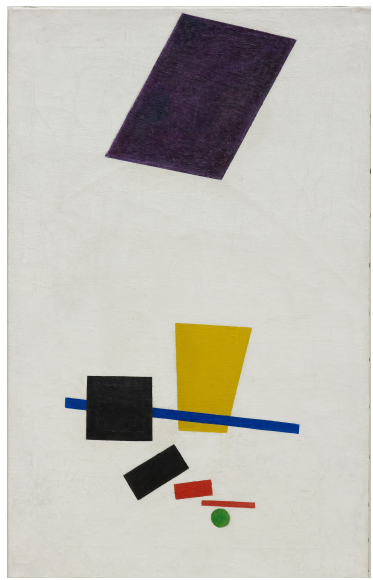

(14) HAFbsVPXhzyLMg sur l'Institut culturel Google / Domaine public

Since that period when abstract art was born, many artists have borrowed mathematical objects and ideas for their work. Let us give just some famous cases. Miró is a clear case of a painter whose paintings contain many geometrical objects. Salvador Dalí's last painting, *The Swallow's Tail* (1983) (Fig. 15), was part of a series inspired by the mathematician René Thom's *catastrophe theory*. Vasarely's *optical illusions* are just based on geometric objects and symmetries. Another famous painter whom we can speak of here is Escher, who was not a mathematician, but who was fascinated by mathematics and whose work uses mathematical objects and operations including impossible objects, explorations of infinity, reflection, symmetry, perspective, truncated and stellated polyhedra, hyperbolic geometry, and tessellations.

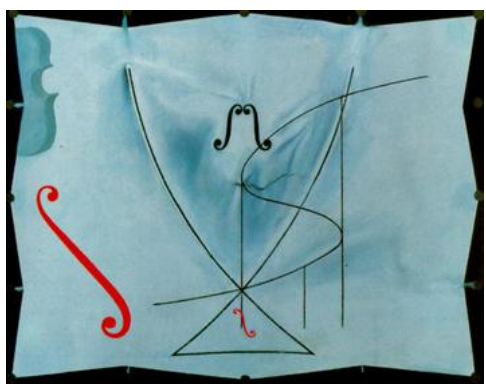

(15) Photo, www.dali-gallery.com / Fair use

Another big interaction between mathematics and art is related to fractals. *Fractal art* (Fig. 16) is a form of algorithmic art created by generating fractal objects applying iterative methods to solve nonlinear equations or polynomial equations. *Fractal art* developed from the mid-1980s onwards, and in recent years its presence has increased with the help of computers because of the calculative capabilities they provide. Surprisingly, Islamic geometric patterns are reminiscent of *fractal art*, as, for instance, in the main dome of *Selimiye Mosque in Edirne* (Fig. 17), Turkey, designed using self-similar patterns. Or in multiple patterns that can be seen in Persian carpets.

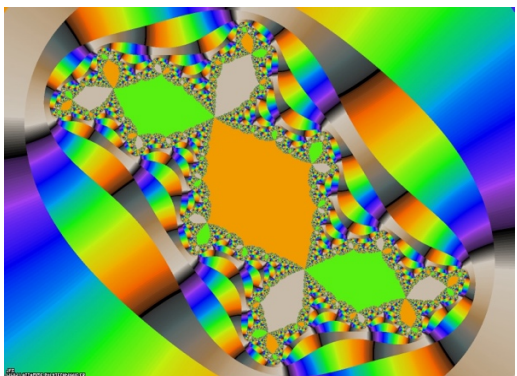 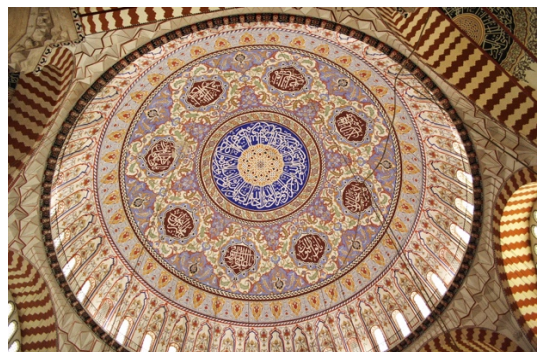

(16) Figure by Jean-François. Colonna (CMAP/Ecole polytechnique)     (17) Photo by Stylommatophora / CC BY 4.0

It is impossible to list or to discuss all the painters that have drawn inspiration from mathematics, and specially from geometry, when painting their artwork. The examples described above give a good idea of the importance of that phenomenon and of the continuity of this practice since the beginning of the 20th Century.

## 5. Mathematics in photography and the movie industry.

### 5.1. Photography of math and math for photography

Very early, photographers were attracted by geometric objects and shapes. Maybe the best example of this is Man Ray who following a suggestion of Max Ernst, photographed in 1934 some 34 objects in the *mathematical model collection* of the *Institut Henri Poincaré* in Paris. These spectacularly beautiful objects were solutions to some algebraic equations. Man Ray's photographs attracted some attention in the artistic world of the time and were quickly featured in art magazines and in exhibitions. After World War II he decided to paint 33 oil paintings based on the original photographs, and he gave to each of them the name of a play by Shakespeare. He called those paintings the 'Shakespearean Equations'. A beautiful movie entitled *Man Ray et les equations shakespeariennes* [10] has been produced recently by the *Institut Henri Poincaré*. It features the mathematical models, Man Ray's photographs and an interview with Man Ray himself. Other interesting documentaries are *Math mon modèle : Man Ray à l'Institut Henri Poincaré - Interview d'Edouard Sebline* [11] and *Mystérieux modèles de Man Ray* [12].

But probably the main relations between photography and mathematics lie elsewhere. In photography mathematics is used in many ways. It is used to understand optical principles, to design lenses, to calculate aperture, lighting, etc. Strikingly, various fields of mathematics are commonly used in photography, like geometry, to calculate angles and distances, algebra and also analysis when dealing with retouching and with blurring. Mathematics is also used to modify photographs, to clean and deblur them, to change the light, the colors, the depth of the field, and other effects.

Mathematics is also key in the design of the universally used JPEG 2000 protocol in image, video and audio compression. The striking improvement from JPEG to JPEG 2000 was possible because of the use of the *wavelet theory*, developed by engineers and mathematicians. The more recent JPEG XS protocol goes still further in the high-quality compression technology, again with the help of wavelets. Applications of JPEG XS include streaming high-quality images for autonomous vehicles, drones, virtual reality, etc.

Finally, mathematics, and especially analysis and partial differential equations, are used in photographs and movies' restoration. The techniques used before 1980 or so to clean and restore images were revolutionized by the discovery of new filters based on sophisticated nonlinear partial differential equations which provided stunning results.

### 5.2. Math in the video and the movie industry today

Believe it or not, mathematics is nowadays ubiquitous in the movie industry. It is used in the design of camera equipment, in the creation of animation movies, in the production of special effects and even in the edition and post-production processes. How can you build a new camera without using a lot of optics, geometry and trigonometry? How calibrate the light or the sound in every scene of a movie without mathematical equations?

Of course, when making a movie, the movie director or the cameramen will not be using mathematics directly, but the material they use has all that incorporated, and the algorithms used in every action have been developed using mathematics and mathematical equations. If this is the case when making a classical movie, with actors in front of the camera in a natural landscape, it becomes more interesting when it comes to animation movies. There you need to design the objects and the virtual actors, and you must *move them*. Doing that requires algorithms based on many different mathematical equations, and again, also on algebra, geometry and trigonometry. But there is a kind of movies that do not belong to the above two categories, and where mathematics is again of immense help. This is when one wants to design special effects without having them for true, because it is too difficult technically, or too dangerous or too expensive. For instance, many of the scenes in the film *Titanic* were not real, in the sense that the filmmakers did not just film real images of water inside the boat, or of people carried away by the water. To give the impression that the images were real, an immense number of computations were run using, for instance, the Navier-Stokes equations to simulate water motion in various different contexts. Most of the water we saw in that movie was not real water, but water motion simulated using mathematics.

In the animation industry there is the need to have good computer models of mechanics, biomechanics, etc., to create beautiful and above all, realistic objects and characters. I remember learning some years ago that one of the most difficult things when creating virtual images is how to represent the motion of hair and textile textures in a realistic way. A very good manner to understand the difficulties around animation and how mathematics helps with describing hair or fabric motion is to watch the documentary *Movie Magic: The Mathematics behind Hollywood's Visual Effects* [13], a lecture given by Eitan Grinspun where he explains how differential geometry and differential equations are all over the place when making an animated movie. Researchers in mathematics, physics and engineering have invested a lot in this direction, and the movie industry has also invested a lot of money funding such research efforts. Note that this kind of research goes much further than making movies, it can be useful, for instance, in medicine, in robotics, in fashion and in virtual reality.

The kind of mathematics used in this area has also evolved very quickly in the last decades. If models based on differential equations were the norm some twenty years ago, nowadays, other methods with mathematics inside, have taken the relay. For instance, *Generative AI*, *NERFS* and the quite recent *Gaussian splatting*, the latest methods in computer graphics.

6. **Mathematics and other arts: Literature, Poetry.**

Literature is an artistic activity where the impact of mathematics is maybe not so obvious. Nevertheless, some of its impact in that field can also be observed, especially in poetry. Apart from the clear numerology, which is behind the rhymes, there are other mathematical objects that are inherent to poetry in many cultures. For instance, prime numbers in the Japanese *haikus*, even numbers in Western poetry and odd numbers in Oriental poetry.

The Bridges conferences association which oversees the annual Bridges conference on mathematical connections in art, music, architecture, and culture organized a special session devoted to poetry [14]. Also, very interesting and original discussions have taken place in a seminar organized at the *Ecole Normale Supérieure* in Paris about *Poésie et Mathématiques. Le fond et la forme* [15].

But well beyond poetry, modern writers have been inspired by mathematics, sometimes in very surprising ways. An example is the *Ouvroir de littérature potentielle (OULIPO),* group of writers seeking to create works using constrained writing techniques. The famous writer Raymond Queneau was one of its founders. He was very attracted to mathematics as a source of inspiration. The group defines the term *littérature potentielle* as "the seeking of new structures and patterns which may be used by writers in any way they enjoy" and they propose to use new methods, often based on mathematical problems. The famous writer Italo Calvino was a member of *OULIPO* until his death. But the *OULIPO* writers are not the only ones imposing constraints when they write a book. For instance, E. Eleanor Catton's Booker Prize-winning *The Luminaries* was written by imposing a precise numerical rule on its chapters, each of which is half the length of the last.

Let us close this section by adhering to what the British mathematician Sarah Hart says in an opinion article published by the *New York Times* [16], "There is a deeper reason we find mathematics at the heart of literature. The universe is full of underlying structure, pattern and regularity, and mathematics is the best tool we have for understanding it — that's why mathematics is often called the language of the universe, and why it is so vital to science. Since we humans are part of the universe, it is only natural that our forms of creative expression, literature among them, will also manifest an inclination for pattern and structure".

7. **Mathematical creation as an art by itself.**

In the above sections we have spoken a lot about how mathematics, or mathematical objects and concepts, are used in art to build, to paint, to design. We have also discussed how mathematics is nowadays used to produce art, like new music, videos, movies, etc. In all those activities mathematics appears as a tool or an inspiration for artists. But what about mathematics itself? One can define art as an activity where beauty is created. But does not mathematics do the same? What is the difference between the creativity of an artist and the creativity of a mathematician developing a new theory or finding a new proof for a beautiful and interesting theorem? In both cases, the actors, artists or mathematicians, use similar spatial reasoning skills and work with patterns.

Mathematics has often been linked to philosophy. Because of logics and because of the necessary rationality and reasoning in both fields. But outside of mathematics few people seem to appreciate the beauty of mathematical creation. Nevertheless, in the previous sections we have often spoken of artists who felt attracted to mathematics, who liked the possibilities that mathematical theories offered to them in their creative activities. They appreciated the mathematical structures and processes that were behind it, or did they only appreciate the geometric shapes and objects, the symmetries, the numbers and operations yielding beautiful proportions?

Mathematicians often speak of a beautiful theory, a beautiful theorem or a beautiful proof. What do they mean by that? I am afraid that there is no general theory about this and certainly no consensus at all. Different mathematicians will appreciate different properties or features when they define beauty in their work of in others' work. It is, like in art, a question of taste. For some, a beautiful proof will be a short, direct and elegant one. For others it will be one which is generated from the interaction of different fields and concepts. For others still, it will be a proof that solves a difficult and important problem. But what is certain is that many mathematicians have an artistic view of themselves and their work. They see themselves as 'creators', and most of them aspire to beauty, even if others will be more interested in usefulness, that is, in work that helps to understand and apprehend nature, physics, biology, health issues, chemistry, or social and economic problems.

It is true that pure mathematics, and more concretely, fields like number theory, topology, knot theory, group theory, etc, are fields that seem to be more akin to beauty in an abstract manner. Symmetry too, and therefore group theory, is a concept that is very linked to beauty, and people working with tilings, periodic or aperiodic, and producing them, feel like real artists producing beautiful designs.

Many mathematicians do their work only in search of interesting and 'beautiful' objects and relations between them. Many of them would even object to do something that is useful, they are only interested in beauty, or what they think beauty is. Some very well-known mathematicians even objected explicitly to any utility of what they did beyond the beauty of it. G.H. Hardy was one of them, even if later what he did has proved to be very useful to solve many applied problems.

Inspiration is something that artists will speak of when explaining their creative work. Mathematicians too. Inspiration is very important when exploring new territory. What is the difference of inspiration for an artist or a mathematician? In my opinion, none that can be clearly explained.

So, many of the characteristics of artistic creation are met by mathematicians when they create their theories, their new theorems and their new visions. In that respect, it is clear to me that mathematical creation is an art!

**Conclusions**

Mathematics is a science, and many people would see it very far away from the arts. But in this article, we have shown, often with concrete examples, that there are strong links between the two. And not only because one sees mathematical objects in artistic paintings, sculptures and buildings, but mainly because mathematical techniques help artists create their art. This was always the case, since the oldest antiquity, when mathematics appeared, among other reasons, to help build big structures or count or measure. Along the centuries, mathematical tools have helped artists to sculpt and paint using principles that allowed their work to be beautiful and realistic. More recently, in the modern era, mathematics has occupied an important place in art, because due to the 20$^{th}$ Century artistic movements tending to simplicity and functionality, mathematical objects have been present in the design of buildings, in the surfaces of paintings and in the shapes of sculptures, and this much more than in the past, because it became explicit. Mathematics abandoned then its role of set of tools, of techniques, to become art. Finally, at the end of the 20$^{th}$ century, the role of mathematics became again that of a set of tools and techniques, since the digital era has armed artists with algorithms and computers with which they can create. But those algorithms are often based on mathematics, and so mathematics is again being used as a basic, often hidden, technology to produce the most advanced artistic creations.

For interested readers, some references not quoted above and which complete this article very well are [17-24].